\documentclass[12pt,reqno]{amsart}
\usepackage{setspace}
\input{amssym.def}
\input{amssym.tex}

\font\bb=msbm10 at 12pt

\newcommand{\bbZ}{\mbox{\bb Z}}

\newtheorem{theorem}{Theorem}
\newtheorem{lemma}[theorem]{Lemma}
\newtheorem{corollary}[theorem]{Corollary}

\theoremstyle{definition}

\def\fix{{\rm Per}}
\def\least{{\rm LPer}}
\def\trace{{\rm{trace}}}

\begin{document}

\doublespacing

\title{A dynamical property unique to
the Lucas sequence}
\author{Yash Puri}
\author{Thomas Ward}
\dedicatory{\rm School of Mathematics,
University of East Anglia,
Norwich,
NR4 7TJ,
U.K.}
\date{\today}
\thanks{The first author gratefully acknowledges
the support of E.P.S.R.C. grant 96001638}
\maketitle

\section{Introduction}

A {\sl dynamical system} is taken here to mean a
homeomorphism
$$
f:X\to X$$
of a compact metric space
$X$ (though the observations here apply equally
well to any bijection on a set).
The number of points with period $n$ under $f$ is
\begin{equation*}
\fix_n(f)=\#\{x\in X\mid f^nx=x\},
\end{equation*} and
the number of points with least period $n$
under $f$ is
\begin{equation*}
\least_n(f)=\#\{x\in X\mid\#\{f^kx\}_{k\in\Bbb Z}=n\}.
\end{equation*}
There are two basic properties that the
resulting sequences $\left(\fix_n(f)\right)$ and
$\left(\least_n(f)\right)$
must satisfy if they are finite. Firstly, the
set of points with period $n$ is the disjoint union of
the sets of points with least period $d$ for each divisor
$d$ of $n$, so
\begin{equation}\label{one}
\fix_n(f)=\sum_{d\vert n}\least_d(f).
\end{equation}
Secondly, if $x$ is a point with least period
$d$, then the $d$ distinct points
$x,f(x),f^2(x),\dots,f^{d-1}(x)$ are all points with
least period $d$, so
\begin{equation}\label{two}
0\le\least_d(f)\equiv 0\mbox{ mod }d.
\end{equation}
Equation (\ref{one}) may be
inverted via the M{\"o}bius inversion formula to give
\begin{equation*}
\least_n(f)=\sum_{d\vert n}\mu(n/d)\fix_d(f),
\end{equation*}
where $\mu(\cdot)$ is the M{\"o}bius function defined
by
\begin{eqnarray*}
\mu(n)
=\left\{
\begin{array}{cl}
1&\mbox{ if }
n=1,\\
0&\mbox{ if $n$ has a squared factor, and}\\
(-1)^r&\mbox{ if $n$ is a product of $r$ distinct primes.}
\end{array}
\right.
\end{eqnarray*}
A short proof of the inversion formula
may be found in \cite[Section 2.6]{wilf-1994}.

Equation (\ref{two}) therefore implies that
\begin{equation}\label{divisibilitycondition}
0\le\sum_{d\vert n}\mu(n/d)\fix_d(f)\equiv 0\mbox{ mod $n$}.
\end{equation}
Indeed,
equation (\ref{divisibilitycondition}) is the
only condition on periodic points in dynamical systems:
define a given sequence of non-negative integers
$\left(U_n\right)$ to be {\sl exactly realizable} if there
is a dynamical system $f:X\to X$ with
$U_n=\fix_n(f)$ for all $n\ge 1$.
Then
$\left(U_n\right)$ is
exactly realizable if and only if
\begin{equation*}
0\le\sum_{d\vert n}\mu(n/d)U_d\equiv 0\mbox{ mod $n$}\mbox{ for all }
n\ge 1,
\end{equation*}
since the realizing map may be constructed as an infinite
permutation using the quantities
$\frac{1}{n}\sum_{d\vert n}\mu(n/d)U_d$ to determine the
number of cycles of length $n$.

Our purpose here is to study sequences of the form
\begin{equation}\label{definefibanoccish}
U_{n+2}=U_{n+1}+U_n,n\ge1,\quad U_1=a,U_2=b,\quad a,b>0
\end{equation}
with the distinguished Fibonacci sequence denoted
$(F_n)$, so
\begin{equation}\label{relation}
U_n=aF_{n-2}+bF_{n-1}\mbox{ for } n\ge 3.
\end{equation}

\begin{theorem}\label{maintheorem}
The sequence $\left(U_n\right)$ defined by
{\rm(\ref{definefibanoccish})} is exactly realizable if
and only if $b=3a$.
\end{theorem}

This result has two parts: the {\sl existence}
of the realizing dynamical system is described first,
which gives many modular corollaries concerning
the Fibonacci numbers. One of these is used in
the {\sl obstruction} part of the result later.
The realizing system is
(essentially) a very familiar and well-known
system, the {\sl golden-mean shift}.

The fact that (up to scalar multiples) the Lucas sequence
$(L_n)$ is the only exactly realizable
sequence satisfying the Fibonacci recurrence relation
to some extent explains the familiar observation that
$(L_n)$ satisfies a great array of
congruences.

Throughout, $n$ will denote a positive integer
and $p,q$ distinct prime numbers.

\section{Existence}

An excellent introduction to the family of
dynamical systems from which the example comes
is the recent book by Lind and Marcus
\cite{lind-marcus}.
Let
$$
X=\left\{{\mathbf x}=
(x_k)\in\{0,1\}^{\Bbb Z}\mid
x_k=1\implies x_{k+1}=0\mbox{ for all }
k\in\bbZ\right\}.
$$
The set $X$ is a compact metric space in a natural
metric (see \cite[Chapter 6]{lind-marcus} for the details).
The set $X$ may also be thought of as the set
of all (infinitely long in both past and future)
itineraries of a journey involving two locations
($0$ and $1$), obeying the rule that from $1$ you must travel
to $0$, and from $0$ you must travel to either $0$ or $1$.
Define the homeomorphism $f:X\to X$ to be the {\sl left shift},
$$
(f({\mathbf x}))_k=x_{k+1}\mbox{ for all }k\in{\bbZ}.
$$
The dynamical system $f:X\to X$ is a simple example
of a {\sl subshift of finite type}.
It is easy to check that the number of points of period
$n$ under this map is given by
\begin{equation}\label{traceformula}
\fix_n(f)=\trace\left(A^n\right)
\end{equation}
where $A=\bmatrix1&1\\1&0\endbmatrix$
(see \cite[Proposition 2.2.12]{lind-marcus};
the $0-1$ entries in the matrix $A$ correspond to the
allowed transitions $0\to 0$ or $1$; $1\to 0$ in the
elements of $X$ thought of as infinitely long
journeys in a graph with vertices $0$ and $1$).

\begin{lemma}\label{existence} If $b=3a$ in
{\rm(\ref{definefibanoccish})}, then the corresponding
sequence is exactly realizable.
\end{lemma}

\begin{proof} A simple induction argument shows that
(\ref{traceformula}) reduces to
$
\fix_n(f)=L_n\mbox{ for }n\ge1,
$
so the case $a=1$ is realized using the golden mean shift itself.
For the general case, let $\bar{X}=X\times B$ where
$B$ is a set with $a$ elements, and define $\bar{f}:\bar{X}\to
\bar{X}$ by $\bar{f}({\mathbf x},y)=(f({\mathbf x}),y)$.
Then
$\fix_n(\bar{f})=a\times\fix_n(f)$ so we are done.
\end{proof}

The relation (\ref{divisibilitycondition}) must as a result
hold for $(L_n)$.

\begin{corollary}\label{congruence}
$\sum_{d\vert n}\mu(n/d)L_d\equiv 0\mbox{ mod }n$
for all $n\ge1$.
\end{corollary}

This has many consequences, a sample of which we list
here. Many of these are of course well-known
(see \cite[Section 2.IV]{ribenboim-records}) or follow easily
from well-known congruences.

\noindent(a) Taking $n=p$ gives
\begin{equation}\label{othercoolthing}
L_p=F_{p-2}+3F_{p-1}\equiv 1\mbox{ mod $p$}.
\end{equation}

\noindent(b) It follows from (a) that
\begin{equation}\label{coolthing}
F_{p-1}\equiv 1\mbox{ mod $p$ }\Leftrightarrow
F_{p-2}\equiv -2\mbox{ mod $p$},
\end{equation}
which will be used below.

\noindent(c) Taking $n=p^k$ gives
\begin{equation}\label{powerthing}
L_{p^k}\equiv
L_{p^{k-1}}\mbox{ mod $p^k$}
\end{equation}
for all primes $p$ and $k\ge 1$.

\noindent(d) Taking $n=pq$ (a product of distinct primes)
gives
$$
L_{pq}+1\equiv
L_p+L_q\mbox{ mod }pq.
$$

\section{Obstruction}

The negative part of Theorem \ref{maintheorem} is
proved as follows. Using some simple modular results on
the Fibonacci numbers, we show that if the sequence
$\left(U_n\right)$ defined by
(\ref{definefibanoccish}) is exactly realizable, then the
property (\ref{divisibilitycondition}) forces the
congruence $b\equiv 3a$ mod $p$ to hold for
infinitely many primes $p$, so
$(U_n)$ is a multiple of $(L_n)$.

\begin{lemma}\label{easything}
For any prime $p$,
$F_{p-1}\equiv 1$ {\rm mod} $p$ if
$p=5m\pm2$.
\end{lemma}

\begin{proof} From Hardy and Wright,
\cite[Theorem 180]{hardy-and-wright},
we have that $F_{p+1}\equiv 0$ mod $p$ if
$p=5m\pm 2$. The identities
$F_{p+1}=2F_{p-1}+F_{p-2}\equiv 0$ mod $p$ and
(\ref{othercoolthing}) imply that
$F_{p-1}\equiv 1$ mod $p$.
\end{proof}

Assume now that the sequence
$\left(U_n\right)$ defined by
(\ref{definefibanoccish}) is exactly realizable.
Applying (\ref{divisibilitycondition})
for $n$ a prime $p$ shows that
$$
U_p-U_1\equiv 0\mbox{ mod }p,
$$
so by (\ref{relation})
\begin{equation*}
aF_{p-2}+bF_{p-1}\equiv a\mbox{ mod $p$}.
\end{equation*}
If $p$ is $2$ or $3$ mod $5$, Lemma \ref{easything}
then implies that
\begin{equation}\label{eight}
\left(F_{p-2}-1\right)a+b\equiv0\mbox{ mod $p$}.
\end{equation}
On the other hand,
for such $p$, (\ref{coolthing}) implies that $F_{p-2}\equiv -2$ mod
$p$, so (\ref{eight}) gives
\begin{equation*}
b\equiv 3a\mbox{ mod $p$}.
\end{equation*}
By Dirichlet's theorem
(or simpler arguments) there are infinitely many primes $p$ with
$p$ equal to $2$ or $3$ mod $5$, so
$b\equiv 3a$ mod $p$ for arbitrarily large values of $p$.
We deduce that $b=3a$, as required.

\section{Remarks}

(a) Notice that the example of the golden mean shift plays
a vital role here. If it were not to hand, exhibiting
a dynamical system with the required properties would
require {\sl proving} Corollary \ref{congruence},
and {\it a priori} we
have no way of guessing or proving this congruence
without using the dynamical system.

(b) The congruence (\ref{othercoolthing}) gives a
different proof that $F_{p-1}\equiv0\mbox{ or }1$ mod $p$
for $p\neq 2,5$. If $F_{p-1}\equiv\alpha$ mod $p$, then
(\ref{othercoolthing}) shows that $F_{p-2}\equiv 1-3\alpha$
mod $p$, so $F_{p}\equiv 1-2\alpha$.
On the other hand, the recurrence relation
gives the well-known equality
$$
F_{p-2}F_{p}=F_{p-1}^2+1,
$$
(since $p$ is odd)
so $1-5\alpha+6\alpha^2\equiv \alpha^2+1$, hence
$5(\alpha^2-\alpha)\equiv 0$ mod $p$. Since
$p\neq 5$, this requires that $\alpha^2\equiv\alpha$ mod $p$
so $\alpha\equiv 0\mbox{ or }1$.

(c) The general picture of conditions on linear
recurrence sequences that allow exact realization is
not clear, but a simple first step in the Fibonacci spirit
is the following question. For each $k\ge 1$ define a recurrence sequence
$(U_n^{(k)})$ by
$$
U_{n+k}^{(k)}=U^{(k)}_{n+k-1}+U^{(k)}_{n+k-2}+\dots+U^{(k)}_{n}
$$
with specified initial conditions $U_j^{(k)}=a_j$ for
$1\le j\le k$.
The subshift of finite type associated to the $0-1$ 
$k\times k$ matrix
$$
A^{(k)}=\bmatrix
1&1&1&\hdots&1&1\\
1&0&0&\hdots&0&0\\
0&1&0&\hdots&0&0\\
&&\ddots\\
0&0&\hdots&1&0&0\\
0&0&\hdots&0&1&0
\endbmatrix
$$
shows that the sequence $(U_n^{(k)})$ is exactly
realizable if
$a_j=2^j-1$ for $1\le j\le k$.
If the sequence is exactly realizable, does
it follow that $a_j=C(2^j-1)$ for $1\le j\le k$ and some
constant $C$? The special case $k=1$ is trivial,
and $k=2$ is the argument above.
Just as in Corollary~\ref{congruence},
an infinite family of congruences follows for
each of these multiple Fibonacci sequences
from the existence of the exact realization.
\bibliographystyle{plain}

\vskip1in

\noindent 1991 {\it Mathematics Subject Classification.}{\ }

\end{document}